\documentclass[10pt]{article}

\usepackage{color}
\usepackage{graphicx}
\usepackage{amsfonts}
\usepackage{hyperref}
\usepackage{natbib}

\def\Eq#1{\label{#1}}

\definecolor{iblue}{RGB}{65,105,225}
\definecolor{ired}{RGB}{220,20,60}
\definecolor{igreen}{RGB}{50,205,50}
\definecolor{ipurple}{RGB}{75,0,130}
\definecolor{iochre}{RGB}{218,165,32}
\definecolor{iteal}{RGB}{51,204,204}
\definecolor{imauve}{RGB}{204,51,153}
\definecolor{RED}{RGB}{255,0,0}

\def\alertb#1{{\color{blue}#1}}
\def\alertr#1{{\color{red}#1}}
\def\alertr#1{{\color{RED}#1}}

\let\a=\alpha    \let\g=\gamma       \let\e=\varepsilon
  \let\h=\eta       \let\l=\lambda
                 \let\p=\pi     
     \let\f=\varphi    
  \let\ps=\psi      \let\o=\omega

\def\V#1{{\bf#1}}
\def\*{\vskip 3mm}
\def\0{\noindent}
\def\be{\begin{equation}}
\def\ee{\end{equation}}
\def\bea{\begin{eqnarray}}
\def\eea{\end{eqnarray}}

\def\EE{{\cal E}}\def\LL{{\mathcal L}}
\def\FF{{\cal F}}

\let\dpr=\partial\let\fra=\frac
\def\ie{{\it i.e.}\ }

\def\otto{\,{\kern-1.truept\leftarrow\kern-5.truept\to\kern-1.truept}\,}

\def\tende#1{\,\vtop{\ialign{##\crcr\rightarrowfill\crcr
 \noalign{\kern-1pt\nointerlineskip} \hskip3.pt${\scriptstyle
 #1}$\hskip3.pt\crcr}}\,}
\def\ie{{\it i.e.\ }}
\def\FF{{\mathcal F}}

\def\ap{{\it a priori}}

\def\({\left(}
\def\){\right)}

\def\iniz{\setcounter{equation}{0}}

\newdimen\xshift \newdimen\xwidth \newdimen\yshift \newdimen\ywidth%
\def\ins#1#2#3{\vbox to0pt{\kern-#2pt\hbox{\kern#1pt #3}\vss}\nointerlineskip}

\def\eqfig#1#2#3#4#5{
\par\xwidth=#1pt \xshift=\hsize \advance\xshift
by-\xwidth \divide\xshift by 2
\yshift=#2pt \divide\yshift by 2%
{\hglue\xshift \vbox to #2pt{\vfil
#3 \includegraphics{#4.eps}
}\hfill\raise\yshift\hbox{#5}}}
\def\Eqfig#1#2#3#4#5#6{
\par\xwidth=#1pt \xshift=\hsize \advance\xshift
by-\xwidth \divide\xshift by 2
\yshift=#2pt \divide\yshift by 2%
{\hglue\xshift \vbox to #2pt{\vfil
#3 \includegraphics{#4.eps}\kern200pt\includegraphics{#5.eps}
}\hfill\raise\yshift\hbox{#6}}}

\def\eqalign#1{\null\,\vcenter{\openup\jot
  \ialign{\strut\hfil$\displaystyle{##}$&$\displaystyle{{}##}$\hfil
      \crcr#1\crcr}}\,}

\def\qedsymbol{$\square$}
\def\qed{\penalty10000\hfill\penalty10000\qedsymbol}

\date{}
\author{\alertb{Giovanni Gallavotti${}^1$ and
Ian Jauslin${}^2$}}

\title{\alertr{\bf A Theorem on Ellipses, an Integrable System and a Theorem
    of Boltzmann}
}

\begin{document}
\maketitle
\kern-8mm
\centerline{${}^1$ INFN-Roma1 \& Universit\`a ``La Sapienza'', email: giovanni.gallavotti@roma1.infn.it}
\centerline{${}^2$ Department of Physics, Princeton University, email: ijauslin@princeton.edu}

\begin{abstract}
We study a mechanical system that was considered by Boltzmann in 1868 in the
context of the derivation of the canonical and microcanonical ensembles. This
system was introduced as an example of ergodic dynamics, which was central to
Boltzmann's derivation. It consists of a single particle in two dimensions,
which is subjected to a gravitational attraction to a fixed center. In
addition, an infinite plane is fixed at some finite distance from the center,
which acts as a hard wall on which the particle collides elastically. Finally,
an extra centrifugal force is added. We will show that, in the absence of this
extra centrifugal force, there are two independent integrals of motion.
Therefore the extra centrifugal force is necessary for Boltzmann's claim of
ergodicity to hold.
\end{abstract}
\*
\0Keywords: {\small Ergodicity, Chaotic hypothesis, Gibbs
distributions, Boltzmann, Integrable systems}
\*

In 1868, \cite{Bo868a} laid the foundations for our modern understanding of the
behavior of many-particle systems by introducing the ``microcanonical
ensemble'' (for more details on this history, see \cite{Ga016}). The principal
idea behind this ensemble is that one can achieve a good understanding of
many-particle systems by focusing not on the dynamics of each individual
particle, but on the statistical properties of the whole. More precisely, the
state of the system becomes a random variable, chosen according to a
probability distribution on phase space, which came to be called the
``microcanonical ensemble''. An important assumption that was made implicitly
by Boltzmann is that the dynamics of the system be ergodic. In this case,
time-averages of the dynamics can be rewritten as averages over phase space,
and the qualitative properties of the dynamics can be formulated as statistical
properties of the microcanonical ensemble.

\begin{figure}[ht]
\hfil\includegraphics[width=6cm]{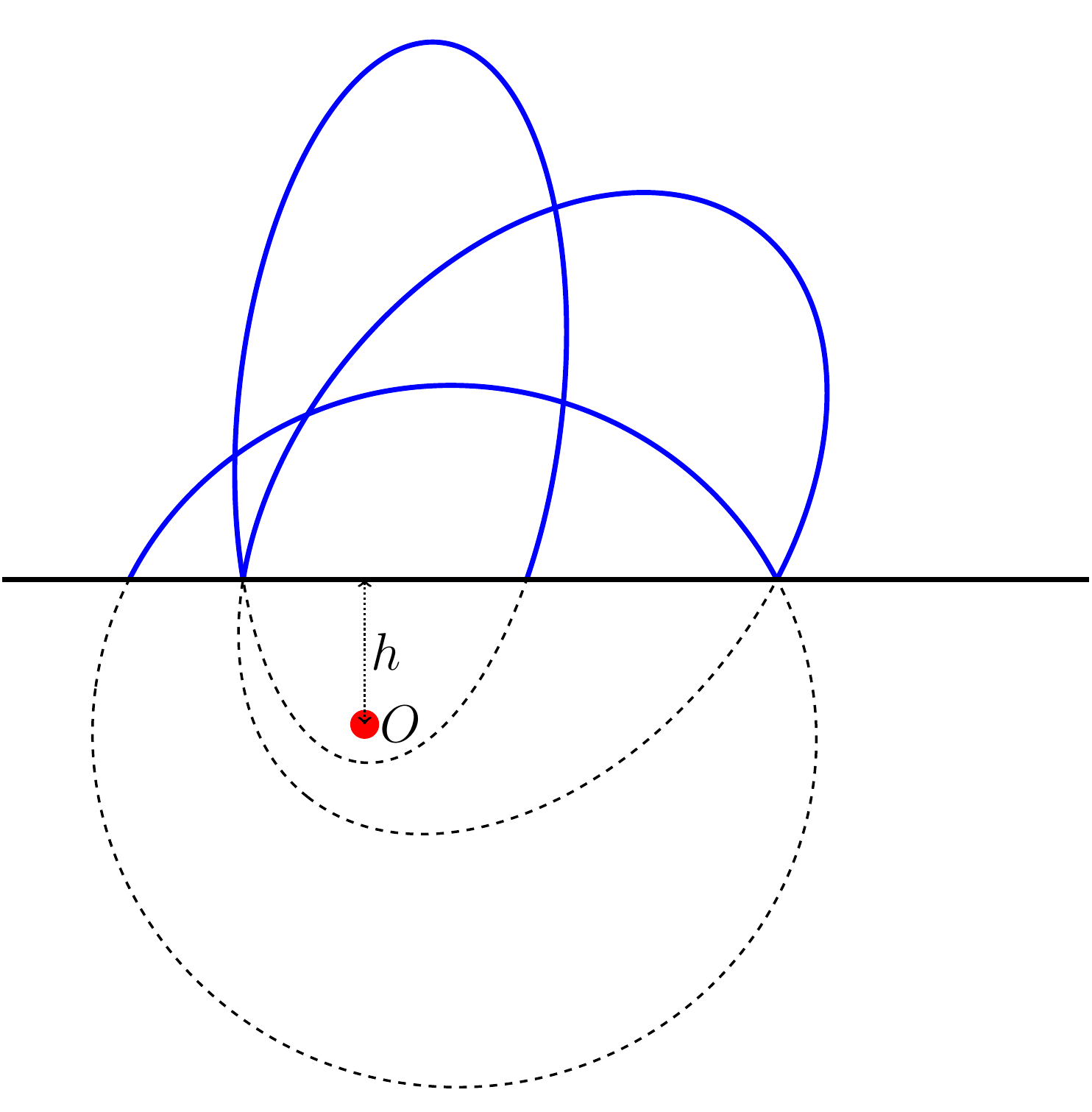}

\caption{A trajectory. The large dot is the attraction center $O$, and the line
  is the hard wall $\LL$. In between collisions, the trajectories are ellipses.
  The ellipses are drawn in full, but the part that is not covered by the
  particle is dashed.
  }
\label{trajectory}
\end{figure}

To support this assumption, Boltzmann presented a mechanical system that
very same year (\cite{Bo868b}) as an example of an ergodic system. This
system consists of a particle in two dimensions that is attracted to a
fixed center via a gravitational potential $-\frac{\alpha}{2r}$. In
addition, he added an extra centrifugal potential $\frac g{2r^2}$. As was
known since at least the times of Kepler, this system is subjected to a
central force, and is therefore integrable. In order to break the
integrability, Boltzmann added an extra ingredient: a rigid infinite planar
wall, located a finite distance away from the center (see figure
\ref{trajectory}). Whenever the particle hits the wall, it undergoes an elastic
collision and is reflected back. Boltzmann's argument was, roughly, that in
the absence of the wall, the dynamics is quasi-periodic, so the particle
should intersect the plane of the wall at points which should fill up a
segment of the wall densely as the dynamics evolves, and concluded that the
region of phase space in which the energy is constant must also be filled
densely. As we will show, this is not the whole story; following a
conjectured integrability for $g=0$, \cite[p.150]{Ga013b},
and first tests
in
\cite[p.225--228]{Ga016}, we have found that, in the absence of the
centrifugal term $g=0$, the dynamics (which has two degrees of freedom) still admits
two constants of motion even in presence of the hard wall.  This
suggests that, if a suitable KAM analysis could be carried out, the system
would not be ergodic for small values of $g$.

\setcounter{section}{0}
\section{Definition of the model and main result}
\label{sec1}
\iniz

Let us now specify the model formally, and state our main result more
precisely. We fix the gravitational center to the origin of the $x,y$-plane and
let $\LL$ be the line $y=h$. The Hamiltonian for the system in between
collisions is
\be H=\frac{p_x^2+p_y^2}2 -\frac{\a}{2r}+\frac{g}{2r^2}
\Eq{e1.1}\ee
where $\a>0,g\ge0,r=\sqrt{x^2+y^2}$ and the particle moves following
Hamilton's equations as long as it stays away from the obstacle $\LL$. When
an encounter with $\LL$ occurs the particle is reflected elastically and
continues on.

\cite{Bo868b}, considered this system on the hyper-surface $A={\V
p}^2-\frac\a r+\frac{g}{r^2}$. The intersection of this hyper-surface with
$y=h$ is the region $\FF_A$ enclosed within the curves
\be \pm\sqrt{(A-\frac{g}{x^2+h^2} +\frac{\a}{\sqrt{x^2+h^2}})},\qquad
x_{min}<x<x_{max}\Eq{e1.2}\ee
with $x_{min}$ and $x_{\max}$ the roots of
$A=\frac{g}{x^2+h^2} -\frac{\a}{\sqrt{x^2+h^2}}$. He argued that
all motions (with few exceptions) would cover densely the surfaces of
constant $A<0$ if $\a,g>0$. 
\bigskip

From now on, unless it is explicitly stated otherwise, we will assume that
$g=0$.
\bigskip

In this case, the motion between collisions takes place at constant
energy $\frac12A$ and constant angular momentum $a$, and traces out an
ellipse. One of the foci of the ellipse is located at the origin, and we
will denote the angle that the aphelion of the ellipse makes with the
$x$-axis by $\theta_0$. Thus, the ellipse is entirely determined by the triplet
$(A,a,\theta_0)$. When a collision occurs, $A$ remains unchanged, but $a$ and
$\theta_0$ change discontinuously to values $(a',\theta_0')=\FF(a,\theta_0)$, and thus
the Kepler ellipse of the trajectory changes. In addition, the semi-major
axis $a_M$ of the ellipse is also fixed to $a_M=-\frac\a{2A}$ (Kepler's
law): so the successive ellipses have the same semi-major axis, while the
eccentricity varies because at each collision the angular momentum changes:
$e^2=1+ \frac{4 A a^2}{\a^2}$. Thus, the motion will take place on arcs of
various ellipses $\EE$, which all share the same focus and the same semi-major
axis, but whose angle and eccentricity changes at each collision.

\begin{figure}[ht]
\hfil\includegraphics[width=100pt]{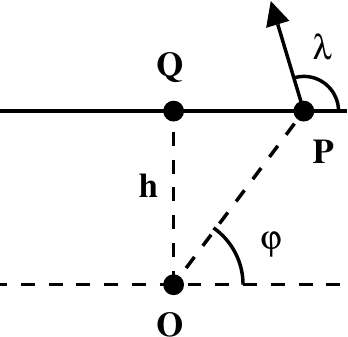}

\caption{The attractive center is $O$, hence it is the focus of the
  ellipse in absence of centrifugal force $g=0$. $Q$ is the projection of
  $O$ on the line $\LL$ and $P$ is a collision point. The arrow
  represents the velocity of the particle after the collision.}
\label{fig1}
\end{figure}

Our main result is that the (canonical) map $(a',\f'_0)=\FF(a,\theta_0)$, which
maps the angular momentum and angle of the aphelion before a collision to their
values after the collision, admits a constant of motion. This follows from the
following geometric lemma about ellipses.
\bigskip

\0{\bf Lemma 1:} {\it Given an ellipse $\mathcal E$ with a focus at $O$ that
intersects $\LL$ at a point $P$. Let $Q$ denote the orthogonal projection of
$O$ onto $\LL$ (see figure \ref{fig1}).  The distance $R_0$ between $Q$ and
the center of $\mathcal E$ depends solely on the semi-major axis $a_M$, the
distance $r$ from $O$ to $P$, and $\cos(2\lambda)$ where $\lambda$ is the angle
between the tangent of the ellipse at $P$ and $\LL$ (to define the direction of
the tangent, we parametrize the ellipse in the counter-clockwise direction):
\be
R_0=\sqrt{\frac14 r^2+\frac14(2a_M-r)^2 +\frac12 r (2a_M-r)\cos(2\l)}
.
\Eq{e1.3}\ee
}
\*

\underline{Proof}: We switch to polar coordinates
$p=(r\cos\f,r\sin\f)$.

\begin{figure}[ht]
  \hfil\includegraphics[width=8cm]{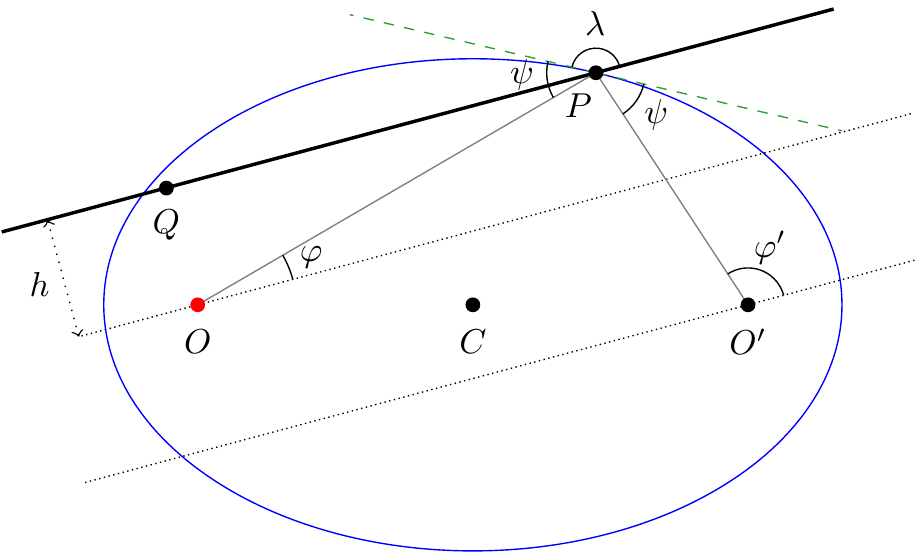}

\caption{An ellipse with foci $O$ and $O'$ and center $C$. The thick line is
$\LL$, which intersects the ellipse at $P$, and $Q$ is the projection of $O$
onto $\LL$. The dashed line is the tangent at $P$. $\lambda$ is the angle
between $\LL$ and the tangent, $\f$ is the polar coordinate, $\f'$ is the angle
between $\LL$ and $\protect\overrightarrow{O'P}$. $\psi$ is the angle between the
tangent and $\protect\overrightarrow{PO}$, which is equal to the angle between the
tangent and $\protect\overrightarrow{PO'}$. $R_0$ is the distance between $Q$ and $C$.}
\label{ellipse}
\end{figure}

Let $O'$ denote the other focus of the ellipse, and $C$ denote its center. The
first step is to compute the vector $\protect\overrightarrow{O'P}$, which in polar
coordinates is
\be
\protect\overrightarrow{O'P}=((2a_M-r)\cos\f',(2a_M-r)\sin\f')\Eq{e1.4}\ee
Let $\psi:=\pi+\f-\lambda$ denote the angle between the tangent of the
ellipse at $P$ and the vector $\protect\overrightarrow{PO}$ (see figure \ref{ellipse}),
and $\psi':=\pi+\f'-\lambda$ denote the angle between the tangent of the
ellipse at $P$ and the vector $\protect\overrightarrow{PO'}$.  

By the focus-to-focus
reflection property of ellipses, we have $\psi'=\pi-\psi$. Thus
$\f'=2\lambda-\pi-\f$ and we find;

\begin{figure}[ht]
  \hfil\includegraphics[width=8cm]{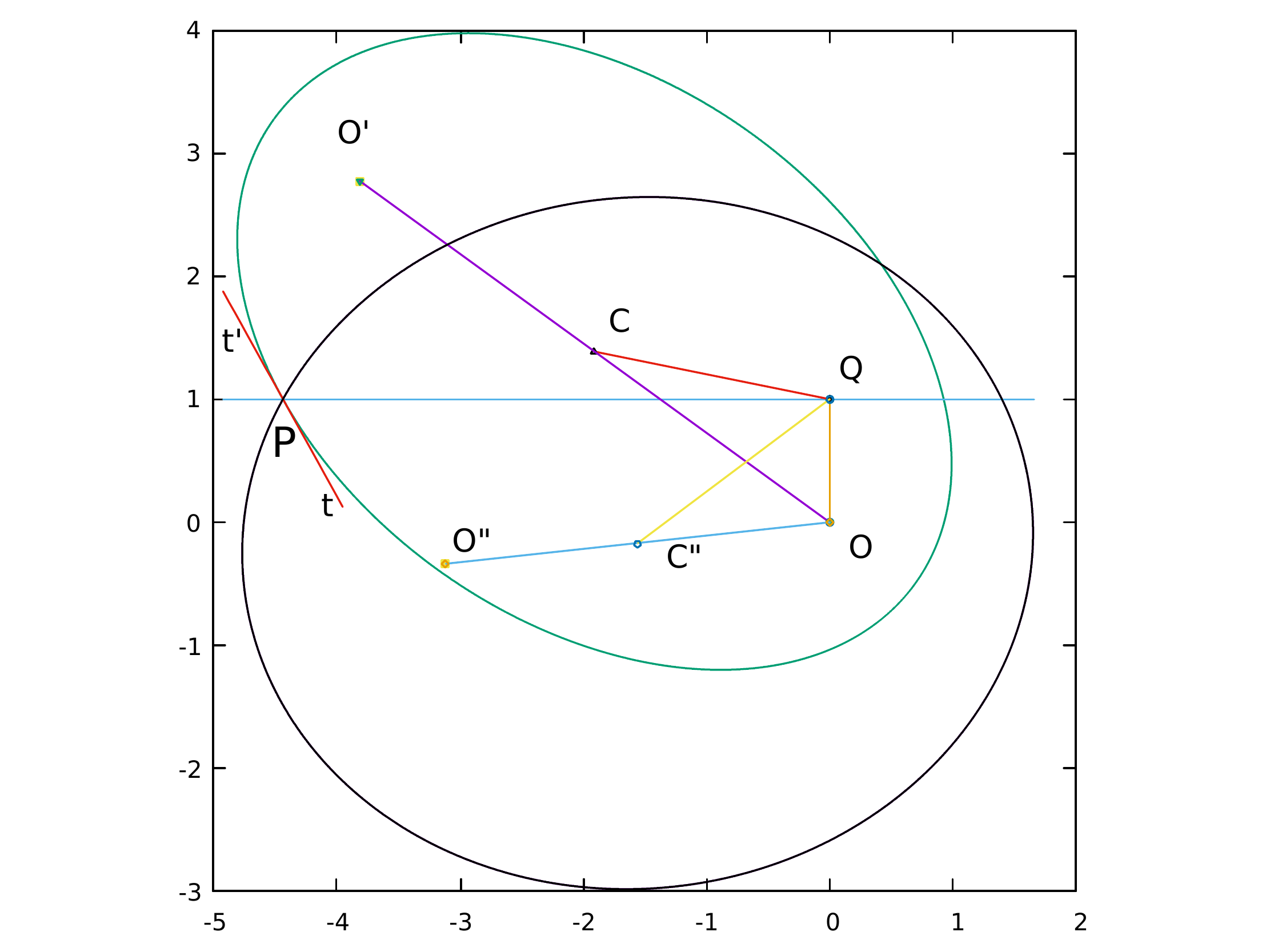}

\caption{Two ellipses, before and after a collision. The collision line $\LL$
is the line at $y=1$, $P$ is the collision point; $Q$ is the projection of $O$
onto $\LL$; the two ellipses $\EE$ and $\EE'$ have a common focus $O$, and
$O,O'$ are the foci of $\EE$, whereas $O,O''$ are the foci of $\EE'$; $C$ and
$C''$ are the centers of $\EE$ and $\EE'$ respectively; the ellipses are drawn
completely although the trajectory is restricted to the parts above $y=h=1$.
The distance from $C''$ to $Q$ is the same as that from $C$ to $Q$.  The upper
ellipse $\EE$ contains the trajectory that starts at the collision point $P$
following the other ellipse $\EE'$ which has undergone reflection.}
\label{fig2}
\end{figure}

\be
  R_0^2=|Q-C|^2
  =
  \frac14\left(r^2+(2a_M-r)^2+2r(2a_M-r)\cos(2\lambda)\right)
  .
\Eq{e1.5}\ee
See figures \ref{ellipse} and \ref{fig2}.\qed

\*
\0{\bf Theorem 1}: {\it The quantity
\be R= a^2+h\a e \sin\theta_0\Eq{e1.6}
\equiv\frac\alpha{2a_M}(h^2+a_M^2-R_0^2)
\Eq{e1.7}\ee
where $e$ is the eccentricity $e=\sqrt{1+\frac{4 A a^2}{\a^2}}$, is a constant
of motion.}
\*

\underline{Proof}:
During a collision, the value of $\l$ changes from $\l$ to $\p-\l$, while
$r$ and $a_M$ stay the same. By lemma 1, this implies that the distance $R_0$
between $Q$ and the center of the ellipse is preserved during a collision.
Furthermore, the position of the center $C$ of the ellipse is given by
$C=a_Me(\cos\theta_0,\sin\theta_0)$
so
\be
  R_0^2=|Q-C|^2=a_M^2e^2-2a_Meh\sin\theta_0+h^2.\Eq{e1.8}
\ee
Furthermore, the angular momentum is equal to
$a^2=\frac12a_M\alpha(1-e^2)$ 
so
\be
  -R_0^2+h^2+a_M^2
  =
  \frac{2a_M}\alpha(a^2+e\alpha h\sin\theta_0)\Eq{e1.9}
\ee
is a conserved quantity.  \qed
\bigskip

\0{\bf Remark:} Some useful inequalities are
\be
\eqalign{
    &r_{max}<{2}{a_M}; \ x_{max}=\sqrt{r_{max}^2-h^2};\
      R_0^2\in ((a_M-r)^2,a_M^2);\cr
      &\frac{\a h^2}{2a_M}\,<\,R\,<\,
      (1+\frac{a_M^2}{h^2}-\Big(\frac{a_M}{h}-
      \frac{r}h\Big)^2)\frac{\a h^2}{2a_M}\cr}
\Eq{e1.10}\ee
hence in the plane $(x,\l)$ the rectangle $(-x_{max},x_{max})\times(0,\p)$
(recall that $x_{max}$ is the largest $x$ accessible at energy $\frac12A$)
is the surface of energy $\frac12A$ and the trajectories are the curves of
constant $R$ inside this rectangle.

\def\SEC{Conjectures on action angle variables}
\section{\SEC}
\label{sec2}
\iniz

In the previous section, we exhibited a constant of motion, which, along
with the conservation of energy, brings the number of independent conserved
quantities to two. In a continuous Hamiltonian system, this would imply the
existence of action-angle variables, which are canonically conjugate to the
position and momentum of the particle, in terms of which the dynamics
reduces to a linear evolution on a torus. In this case, the collision
with the wall introduces some discreteness into the problem, and the
existence of the action angle variables is not guaranteed by standard
theorems. Indeed, in the presence of the collisions, we no longer have a
Hamiltonian system, but rather a discrete symplectic map (or a
non-differentiable Hamiltonian), which describes the change in the state of
the particle during a collision. In this section, we present some
conjectures pertaining to the existence of action angle variables for this
problem.
\bigskip

The first step is to change to variables which are action-angle variables for
the motion in between collision. We choose the {\it Delaunay} variables, whose
angles are the argument of the aphelion $\theta_0$ defined above, the {\it mean
anomaly} $M$, and whose actions are the angular momentum $a$, and another
momentum usually denoted by $L$ and related to the semi-major axis $a_M$ and
to the energy $E=\frac12 A$:
\be L:=-\sqrt{\frac{\alpha}2a_M},\quad a_M:=-\frac\alpha{2A}
,\quad
A:=p^2+\frac{a^2}{r^2}-\frac\alpha{r}\equiv-\frac{\alpha^2}{4L^2}
\Eq{e2.1}\ee
It is well known that this change of variables is canonical. In between
collisions, the dynamics of the particle in the variables
$(M,\theta_0;L,a)$ is, simply,
\be
  \dot M=\frac{\alpha^2}{4L^3}
  ,\quad
  \dot\theta_0=0
  ,\quad
  \dot L=0
  ,\quad
  \dot a=0
  .\Eq{e2.2}
\ee
These variables are thus action-angle variables in between collisions, but when
a collision occurs, $\theta_0$ and $a$ will change.
\bigskip

The following conjecture states that there exists an action-angle variable
during the collisions.
\bigskip

\0{\bf Conjecture 1:} {\it There exists a variable $\gamma$ and an integer
  $k$ such that, every $k$ collisions, the change in $\gamma$ is
  \begin{equation}
    \gamma'=\gamma+\o(L,R)\Eq{e2.3}
  \end{equation}
  in which case $\gamma$ is an angle that rotates on a circle of radius depending
  on $L,R$. The function $\o(L,R)$ has a non zero derivative with respect to
  $R$ at constant $L$, \ie the motion on the energy surface is quasi periodic
  and anisochronous.}
\*

We will now sketch a construction of this variable $\gamma$, which we obtain
using a generating function $F(L,R,M,\theta_0)$.
\bigskip

First of all, by theorem 1, the angular momentum $a(\theta_0)$ is a solution of
\begin{equation}
  a^2=R-h\a\sin\theta_0\sqrt{1-\fra{a^2}{L^2}}\Eq{e2.4}
\end{equation}
that is, if $\e=\pm$,
\be
a^2=R-\frac{h^2\alpha^2}{2L^2}\sin^2\theta_0+
\e\sqrt{\frac{h^4\alpha^2}{4L^4}\sin^4\theta_0+h^2\a^2\sin^2\theta_0-
\frac{R\alpha^2h^2}{L^2}\sin^2\theta_0}\Eq{e2.5}\ee
and $a=\h\sqrt{a^2}$, so that there may be four possibilities for the value of
$a$ denoted $a=a_{\e,\h}(\theta_0,R,L)$ with $\e=\pm,\h=\pm$. The choice of the
signs $\e=\pm1$, and $\h$ must be examined carefully.
\bigskip

We then define the generating function
\be F(L,R,M,\theta_0)=LM+\int_0^{\theta_0} a(L,R,\ps)d\ps
\Eq{e2.6}\ee
which yields the following canonical transformation:

\be\eqalign{
  \g=&\dpr_R\int^{\theta_0}_0 a_{\e,\h}(L,R,\ps)\,d\ps
  \cr
  M'=&M+\dpr_L\int^{\theta_0}_0 a_{\e,\h}(L,R,\ps)\,d\ps
  \cr
}\Eq{e2.7}\ee
It is natural, if Boltzmann's system is integrable (at $g=0$), that the new
variables are its action angle variables and $M',\gamma$ rotate uniformly in spite of
the collisions.
\bigskip

\begin{figure}
  \hfil\includegraphics[width=8cm]{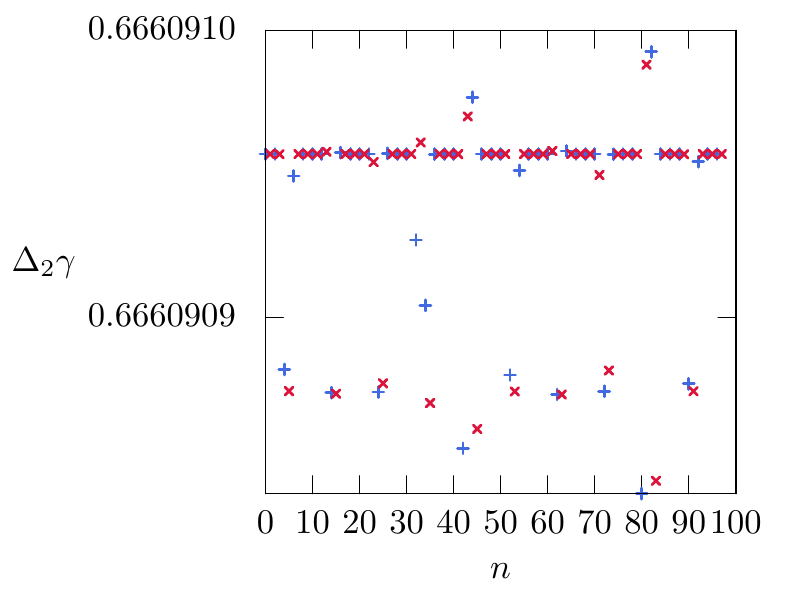}
  \caption{A plot of the increment in $\gamma$ between the $n$-th and the $n+2$-nd collision as a function of $n$. The blue `+' signs correspond to even $n$, and the red `$\times$' to odd $n$. The variation of $\Delta_2\gamma$ is as small as 1 part per million, thus supporting conjecture 2.}
  \label{action_angle}
\end{figure}

However, in this case, the signs $\e$ and $\h$ may change from one collision to
the next, complicating the situation. A careful numerical study of the system
has led us to the following conjecture (see figure \ref{action_angle}).
\bigskip

\0{\bf Conjecture 2:} {\it If $R>h\a$ (which is the case in which the
  circle, of radius $R_0$, of the centers encloses the focus $O$), when the
  motion collides for the $n$-th time, the angular momentum is proportional
  to $(-1)^n$, and, thus, $\epsilon=(-1)^n$. The sign $\eta$ is fixed to
  $+$. The increment $\Delta_2\gamma$ in $\gamma$ between the $n$-th and
  the $n+2$-th collision is independent of $n$.  } \*
\bigskip

\0{\it Remark:} The change of variables over the variables $a,\theta_0$ to
$R,\g$ at fixed $L$ is {\it remarkably} essentially the same as the one
(\ap\ unrelated) to find action-angle variable for the auxiliary
Hamiltonian $R=R(a,\theta_0)$. This might remain true even when $R<h\a$:
interpretable as a kind of auxiliary pendulum motion.  \*

At the time of publication, it has been
brought to our attention that G. Felder has proved that the orbits are all
either periodic or quasi-periodic, which would be implied from conjecture
1.

\section{Conclusion and outlook}
In this brief note, we have shown that the system considered by Boltzmann in
1868, in the case $g=0$, admits two independent constants of motion.
This indicates that it should be possible to compute action
angle variables for this system, which is not entirely trivial because of the
discontinuous nature of the collision process. If such a construction could be
brought to its conclusion, then it would show that the trajectories are either
periodic or quasi-periodic, a fact which is consistent with the numerical
simulations we have run.

This is not a contradiction of Boltzmann's claim that this model is ergodic,
since Boltzmann considered the model at $g\neq 0$. However, we expect that a
KAM-type argument can be set up for this model, to show that the system cannot
be ergodic, even if $g>0$, provided $g$ is sufficiently small. However it
may still have invariant regions of positive volume where the motion is ergodic. 

\*

\0{\bf Acknowledgements}: The authors thank G. Felder for giving us the impetus
to write this note up in its current form, and to publish it. I.J. gratefully
acknowledges support from NSF grants 31128155 and 1802170.

\bibliographystyle{plainnat}

\begin{thebibliography}{3}
\providecommand{\natexlab}[1]{#1}
\providecommand{\url}[1]{\texttt{#1}}
\expandafter\ifx\csname urlstyle\endcsname\relax
  \providecommand{\doi}[1]{doi: #1}\else
  \providecommand{\doi}{doi: \begingroup \urlstyle{rm}\Url}\fi

\bibitem[Boltzmann(1868a)]{Bo868a}
L.~Boltzmann.
\newblock Studien \"uber das Gleichgewicht der lebendingen Kraft zwischen bewegten materiellen Punkten.
\newblock \emph{Wiener Berichte}, {\bf 58}, 517--560, (49--96), 1868.

\bibitem[Boltzmann(1868b)]{Bo868b}
L.~Boltzmann.
\newblock L{\"o}sung eines mechanischen problems.
\newblock \emph{Wiener Berichte}, {\bf 58}, (W.A.,\#6):\penalty0 1035--1044,
  (97--105), 1868.

\bibitem[Gallavotti(2014)]{Ga013b}
G.~Gallavotti.
\newblock \emph{Nonequilibrium and irreversibility}.
\newblock Theoretical and Mathematical Physics. Springer-Verlag, 2014.

\bibitem[Gallavotti(2016)]{Ga016}
G.~Gallavotti.
\newblock Ergodicity: a historical perspective. equilibrium and nonequilibrium.
\newblock \emph{European Physics Journal H}, {\bf 41}, 181--259, 2016.
\newblock \doi{DOI: 10.1140/epjh/e2016-70030-8}.

\end{thebibliography}

\*\*
\end{document}